\documentclass[12pt]{article}
\usepackage{amsmath,amssymb,amsthm,graphics,color}
\newcommand{\Z}{{\bf Z}}
\newcommand{\Q}{{\bf Q}}
\newcommand{\N}{{\bf N}}
\newcommand{\C}{{\bf C}}
\begin{document}

\begin{center}
{\bf Tiling Lattices with Sublattices, I}
\end{center}

\begin{center}
David Feldman, University of New Hampshire \\
James Propp, University of Massachusetts Lowell \\
Sinai Robins, Nanyang Technological University
\end{center}

\noindent
{\small
{\sc Abstract}:
Call a coset $C$ of a subgroup of $\Z^d$ a {\em Cartesian coset}
if $C$ equals the Cartesian product of $d$ arithmetic progressions.
Generalizing Mirsky-Newman, we show that
a non-trivial disjoint family of Cartesian cosets with union
$\Z^d$ always contains two cosets that differ only by translation.
Where Mirsky-Newman's proof (for $d=1$) uses complex analysis, we
employ Fourier techniques.  Relaxing the Cartesian requirement, 
for $d>2$ we provide examples where $\Z^d$ occurs as the disjoint union 
of four cosets of distinct subgroups (with one not Cartesian).
Whether one can do the same for $d=2$ remains open.
}

\bigskip

In 1950, Paul Erd\H{o}s  conjectured that if a system of
$k$ arithmetic progressions, with distinct differences $n_i$,
covers the natural numbers {$\N$} (or the integers $\Z$), then
$$\sum_{i=1}^{k}\frac{1}{n_i}>1\ .$$
Erd\H{o}s credits Leon Mirsky and Donald Newman with the first proof
(and also cites Harold Davenport and Richard Rado for independently
finding the same proof just a bit later).  
Mirsky and Newman did not publish at the time,
but Erd\H{o}s credits them in a 1952 paper.  For Newman's own exposition,
see his analytic number theory monograph.  For more history, see
Soifer.  The Mirsky-Newman proof
blends formal power series and complex analysis, and now serves as
a standard example showing how such interactions solve combinatorial
problems.

Nowadays one usually hears the Mirsky-Newman stated
in the contrapositive: a non-trivial finite disjoint family of arithmetic
progressions with union $\Z$ (so that $\sum_{i=1}^{k}\frac{1}{n_i}=1$)
contains two progressions with the same common difference.
Henceforth we refer to any nontrivial finite disjoint family ${\cal F}$
of sets as a {\em tiling} of its union $\cup {\cal F}$ and speak of the sets
in ${\cal F}$ as {\em tiles}.  Generally speaking, one may ask for groups
$G$ that possess no tiling by cosets of distinct subgroups.
Mirsky-Newman says $G=\Z$ constitutes one such group, but examples we
give below show that the groups $\Z^d$ for $d>2$ do not have this property 
(the case $d=2$ remains open).

We now frame our positive result extending Mirsky-Newman
to $\Z^d$, and obtained by restricting the permissible tiles to cosets of
certain special subgroups.\footnote{Where we strengthen the hypothesis,
Sun's theorem generalizes Mirsky-Newman to $\Z^d$ by weakening the conclusion:
every tiling of $\Z^d$ contains two tiles,
cosets respectively of subgroups $L$ and $L'$ 
giving rise to isomorphic quotients of $\Z^d$.
Whether or not we can require that $L$ and $L'$ 
differ merely by a rotation (in $O(d)$) we do not know.}
We call any coset $T={\bf v}+L$ of a subgroup
$L$ of $\Z^d$ ($d \geq 1$) a {\em Cartesian coset}
if $L$ has the form $a_1 \Z \times \dots \times a_d \Z$
for positive integers $a_1,\dots,a_d$.
Thus a Cartesian coset in $\Z^d$ equals the Cartesian product of $d$
arithmetic progressions.  (Note that $L$ must have full rank.)

\smallskip
\noindent
{\sc Theorem}:
$\Z^d$ admits no tiling using only Cartesian cosets
no two of which are translates of one another.

\smallskip
Ahead of the proof, we gather a few facts of Fourier analysis.
For $\C$, the complex numbers,
every $\C$-valued function on a finite abelian group equals
a unique finite $\C$-linear combination of characters (i.e.\ homomorphisms
to $\C^*$).  For $L$ a full-rank subgroup of $\Z^d$, we
may regard $L$-periodic functions on $\Z^d$ as functions
on $\Z^d/L$ which thus expand as $\C$-linear combinations of characters.
These characters on $\Z^d/L$ pull back to $\Z^d$ as exponential functions
${\bf x} \mapsto \exp(2\pi i {\bf k} \cdot {\bf x})$ 
with ${\bf k}$ such that ${\bf k} \cdot L \subseteq \Z^d$.  
Adding an element of $\Z^d$ to
$\bf k$ will not affect the exponential, so we standardize the components
of ${\bf k}$ to lie in $[0,1)$.  Given $L$, we write $\widehat{L}$
for the associated finite set of standardized vectors ${\bf k}$.
E.g., if $L = 2\Z \times 3\Z$ we have 
$\widehat{L} = \{(0,0),(0,1/3),(0,2/3),(1/2,0),(1/2,1/3),(1/2,2/3)\}$.
% Have I got this right?
Note that the cardinality of $\widehat{L}$ is precisely $[\Z^d:L]$.
% Have I got this right?

Given an $L$-periodic function 
$f: {\bf x} \mapsto \sum_{{\bf k} \in \widehat{L}} 
c_{\bf k} \exp (2 \pi i {\bf k} \cdot {\bf x})$,
write $\widehat{f}({\bf k})=c_{\bf k}$.
% We may recover each $c_{\bf k}$ from $f$
% by an appropriate inverse Fourier transform:
One can write $\widehat{f}({\bf k})$
as $\frac{1}{[\Z^d:L]} \sum_{{\bf x}}
f({\bf x}) \exp(-2 \pi i {\bf k} \cdot {\bf x})$
where ${\bf x}$ varies over a set of coset representatives.
% , or as the limit of
% $\frac{1}{(2N+1)^d} \sum_{{\bf x} \in B_N}
% f({\bf x}) \exp(-{\bf k} \cdot {\bf x})$
% as $N \rightarrow \infty$, where $B_N$ is the box $[-N,N]^d$.
If $f$  is $L$-periodic, then $\widehat{f}({\bf k})$ vanishes
for ${\bf k} \not\in \widehat{L}$.

As usual, for each $S\subseteq Z^d$ we have an indicator function $\chi_S$, 
with $\chi_{S}(v)=1$ if $v\in S$, and $=0$ otherwise.

By symmetry, the sum of all characters on $Z^d/L$ vanishes everywhere 
except the identity, where the value equals
the cardinality of $Z^d/L$; thus summing the corresponding exponentials and
dividing by the index $[\Z^d\!:\!L]$ of $L$ in $\Z^d$ yields $\chi_L$.  
So $\widehat{\chi_L}= 1/[\Z^d\!:\!L]$ on $\widehat{L}$, $0$ elsewhere.

$\chi_{{\bf v}+L}({\bf x})=\chi_{L}({\bf x}-{\bf v})$,
so $\widehat{\chi_{{\bf v}+L}}({\bf k})= (1/[\Z^d\!:\!L])
\exp (-2 \pi i {\bf k} \cdot {\bf v})$ on $\widehat{L}$, 0 elsewhere.

\smallskip
\noindent
{\sc Proof}:
Suppose that $\Z^d=\cup T_j$ gives
a tiling by Cartesian cosets $T_j ={\bf v}_j+L_j$.
Tiling implies $\sum_j \chi_{T_j}=\chi_{\Z^d}$
and thus $\sum_j \widehat{\chi_{T_j}}=\widehat{\chi_{\Z^d}}$.
For $L_j = a_1 \Z \times \dots \times a_d \Z$,
${\bf k}_j := (1/a_1,\dots,1/a_d)\in \widehat{L_j}$.
Take $L_m$ with $a_1 \cdots a_d$ maximal.
Then $\widehat{\chi_{T_m}}({\bf k}_m) \neq 0$ even though
$\widehat{\chi_{\Z^d}}({\bf k}_m)=0$.
So $\widehat{\chi_{T_j}}({\bf k}_m) \neq 0$ for some $j \neq m$.
Our choice of $L_m$ implies $L_j = L_m$.
\hfill $\Box$

\smallskip
\noindent
{\sc Note:}
We know no previous Fourier-theoretic proof of Mirsky-Newman.
However, Henry Cohn has pointed out that in the case $d=1$,
the Fourier coefficients in our proof are precisely
the residues at the poles of the terms in
the partial fraction decomposition that occurs
in the Mirsky-Newman power series proof.
Indeed, the link between Fourier decomposition
and partial fraction decomposition
is discussed in chapter 7 of the book of Beck and Robins.
This link extends to higher dimensions,
and we intend to write a follow-up article
that re-proves the main result of the present paper
using $d$-variable generating functions.
For a preview, see the Feldman-Propp-Robins preprint
listed in the References.

The four cosets
$$(0,1,0)+2\Z\times 2\Z\times \Z,$$
$$(0,0,1)+\Z\times 2\Z\times 2\Z,$$
$$(1,0,0)+2\Z\times \Z\times 2\Z,\ \rm{and}$$
$$(2\Z\times 2\Z\times 2\Z) \cup ((1,1,1)+2\Z\times 2\Z\times 2\Z)$$ 
tile $\Z^3$, and no two of them are translates of one another;
this shows that for $d=3$,
the Theorem becomes false if the Cartesian requirement is dropped.
Furthmore, taking this example $\times \Z^{d-3}$ gives
a $d$-dimensional counterexample for every $d > 3$.

\noindent
{\small
{\sc Acknowledgments}: Conversations with Henry Cohn and suggestions from 
the referee were helpful in the preparation of this article for publication.
The second author is supported by a grant from the National Science Foundation,
and the third author is supported in part by Nanyang Technological University's
SUG grant M5811053.
}

\bigskip

\noindent
{\sc References}

\noindent
M. Beck and S. Robins,
Computing the continuous discretely: integer-point enumeration in polyhedra,
Springer, 2007.

\noindent
P. Erd\H{o}s, Egy kongruenciarendszerekr\H{o}l sz\'{o}l\'{o} probl\'{e}m\'{a}r\'{o}l
(On a problem concerning covering systems). Mat. Lapok 3 (1952), 
122--128 (Hungarian).

\noindent
D. Feldman, J. Propp, and S. Robins,
Tiling Lattices with Sublattices, II,
arXiv: {\tt 1006.0472}.

\smallskip
\noindent
D.J. Newman, Analytic Number Theory, Springer 2000.

\smallskip

\noindent
A. Soifer, The Mathematical Coloring Book, Springer 2008.

\smallskip

\noindent
Z. Sun, On the Herzog-Sch\"onheim conjecture
for uniform covers of groups,
J.\ Algebra 273 (2004), 153--175.

\end{document}